\documentclass[12pt]{article}

\usepackage{amsmath,amssymb,bm,amsthm,mathrsfs,amscd}
\usepackage{color,graphicx}
\newtheorem{Def}{Def}[section]
\newtheorem{Defi}[Def]{Definition}
\newtheorem{Them}[Def]{Theorem}
\newtheorem{Lem}[Def]{Lemma}
\newtheorem{Cor}[Def]{Corollary}
\newtheorem{Prop}[Def]{Proposition}
\newtheorem{Examp}[Def]{Example}

\newtheorem{Claim}[Def]{Claim}

\newtheorem{Prob}[Def]{Problem}

\numberwithin{equation}{section}

\newcommand{\Sym}{\operatorname{Sym}}

\title{A characterization of high order freeness for product arrangements and answers to Holm's questions}
\author{Takuro Abe\footnote{Institute of Mathematics for Industry, Kyushu University, Fukuoka 819-0395, Japan. Email: abe@imi.kyushu-u.ac.jp}\quad
and\quad
Norihiro Nakashima\footnote{Department of Mathematics, Tokyo Denki University, Tokyo 120-8551, Japan. Email: nakashima@mail.dendai.ac.jp}}

\date{}
\begin{document}

\maketitle

\begin{abstract}
An $m$-free hyperplane arrangement is a generalization of a free arrangement.
Holm asked the following two questions:
(1)Does $m$-free imply $(m+1)$-free for any arrangement?
(2)Are all arrangements $m$-free for $m$ large enough?
In this paper, we characterize $m$-freeness for product arrangements,
while we prove that all localizations of an $m$-free arrangement are $m$-free.
From these results, we give answers to Holm's questions.

\vspace{2mm}

\end{abstract}

\section{Main results}
Let $\mathbb{K}$ be a field of characteristic zero, and let
$V$ be an $\ell$-dimensional vector space over $\mathbb{K}$.
A {\bf (central hyperplane) arrangement} $\mathscr{A}=(\mathscr{A},V)$
is a finite collection of hyperplanes in $V$ which contain the origin.
We call $\mathscr{A}$ an $\ell$-arrangement
when we emphasize the dimension of $V$.
For any hyperplane $H\in\mathscr{A}$,
there exists a linear form $\alpha_{H}$
in the dual space $V^{\ast}$ of $V$
such that $\{\alpha_{H}=0\}=H$.
We call $Q=Q(\mathscr{A})=\prod_{H\in\mathscr{A}}\alpha_{H}$
a defining polynomial of $\mathscr{A}$.

Let $\{x_{1},\dots,x_{\ell}\}$ be a basis for the dual space $V^{\ast}$
over $\mathbb{K}$, and let $S=\Sym(V^{\ast})=\mathbb{K}[x_{1},\dots,x_{\ell}]$.
We consider $x_i$ and $\partial_{i}=\frac{\partial}{\partial x_{i}}$
as elements of the endomorphism ring ${\rm End}_{\mathbb{K}}(S)$:
$x_i(f)=x_i f$ and $\partial_{i}(f)=\frac{\partial f}{\partial x_{i}}$
$(f\in S)$.
Let $\mathbb{N}=\{0,1,2,\dots\}$ be the set of nonnegative integers,
and we use multi-index notations:
for $\bm{a}=(a_1,\dots,a_{\ell})\in\mathbb{N}^{\ell}$,
\begin{align}
|\bm{a}|=a_1 +\cdots+a_{\ell},\ 
\bm{a}!=a_1 !\cdots a_{\ell}!,\ 
x^{\bm{a}}=x_1^{a_1}\cdots x_{\ell}^{a_{\ell}}\ {\rm and}\ 
\partial^{\bm{a}}=\partial_{1}^{a_{1}}\cdots\partial_{\ell}^{a_{\ell}}.
\end{align}
Then $D^{(m)}(S)$ is an $S$-submodule\footnote{The ring ${\rm End}_{\mathbb{K}}(S)$ also has an $S$-module structure by multiplying polynomials from the left: $f\cdot\theta=f\theta$ ($f\in S,\ \theta\in{\rm End}_{\mathbb{K}}(S)$).}
of ${\rm End}_{\mathbb{K}}(S)$ defined by
$D^{(m)}(S)=\bigoplus_{|\bm{a}|=m}S\partial^{\bm{a}}$
for $m\geq 1$, and $D^{(0)}(S)=S$.
\begin{Defi}
An $S$-submodule $D^{(m)}(\mathscr{A})$ of $D^{(m)}(S)$ is defined by
\begin{align}
D^{(m)}(\mathscr{A})=\left\{\theta\in D^{(m)}(S)\mid
\theta(QS)\subseteq QS\right\}.
\end{align}
We call $D^{(m)}(\mathscr{A})$
the {\bf module of $m$th order $\mathscr{A}$-differential operators}.
We say that $\mathscr{A}$ is {\bf $m$-free}
if $D^{(m)}(\mathscr{A})$ is a free $S$-module.
\end{Defi}
Since $QD^{(m)}(S)\subseteq D^{(m)}(\mathscr{A})\subseteq D^{(m)}(S)$,
the rank of $D^{(m)}(\mathscr{A})$ is $s=s_{m}(\ell)=\binom{\ell+m-1}{m}$
if $\mathscr{A}$ is $m$-free.
When $m=1$, $D^{(1)}(\mathscr{A})$ is called the module of 
$\mathscr{A}$-derivations which is investigated, relating with
geometries and combinatorics of hyperplane arrangements.
We say that
$\theta=\sum_{|\bm{a}|=m}f_{\bm{a}}\partial^{\bm{a}}\in D^{(m)}(S)$
is {\bf homogeneous of degree $i$} and write $\deg(\theta)=i$,
if $f_{\bm{a}}$ is zero or homogeneous of degree $i$ for each $\bm{a}$.
Then $D^{(m)}(\mathscr{A})$ becomes a graded $S$-submodule of $D^{(m)}(S)$
similar to the derivation module:
$D^{(m)}(\mathscr{A})=\bigoplus_{i\in\mathbb{Z}}D^{(m)}(\mathscr{A})_i$,
where
$D^{(m)}(\mathscr{A})_i=\{\theta\in D^{(m)}(\mathscr{A})\mid\deg(\theta)=i\}$
for $i\geq 0$ and $D^{(m)}(\mathscr{A})_i=\{0\}$ for $i<0$.
If $m\geq 1$ and if $\mathscr{A}$ is $m$-free
with a homogeneous basis $\{\theta_1,\dots,\theta_s\}$,
we define {\bf $m$-exponents} by the multi-set
$\exp_m(\mathscr{A})=\{\deg(\theta_1),\dots,\deg(\theta_s)\}$.
We also define $\exp_0(\mathscr{A})=\{0\}$ for any arrangement $\mathscr{A}$.
Then $m$-exponents depend only on $\mathscr{A}$.

Let $\mathscr{D}(S)$ be the $S$-subalgebra of ${\rm End}_{\mathbb{K}}(S)$
generated by the derivations $\partial_{1},\dots,\partial_{\ell}$,
i.e., $\mathscr{D}(S)$ is the Weyl algebra.
For an ideal $I$ of $S$, let
$\mathscr{D}(I)=\{\theta\in\mathscr{D}(S)\mid \theta(I)\subseteq I\}$.
Then $\mathscr{D}(S/I)=\mathscr{D}(I)/I\mathscr{D}(S)$ is
the set of differentials of $S/I$ (see \cite{Mac-Rob}),
and is called the ring of differential operators of $S/I$.
When $S/I$ is a regular ring,
$\mathscr{D}(S/I)$ has a similar structure to the Weyl algebra,
that is, $\mathscr{D}(S/I)$ is the $(S/I)$-algebra
generated by derivations of $S/I$ (see \cite[Corollary 15.5.6]{Mac-Rob}).
However, as Traves in \cite{Traves-nakaiconj} proved,
$\mathscr{D}(S/I)$ is not generated by derivations
when $S/I$ is a reduced algebra (including the case when $I=QS$).
To observe generators and further structures of $\mathscr{D}(S/QS)$,
Holm in \cite{Holm-phd} proved that
$\mathscr{D}(S/QS)=\sum_{m\geq 0}D^{(m)}(\mathscr{A})/QD^{(m)}(S)$.
A free basis for $D^{(m)}(\mathscr{A})$
is useful to study the structures of $\mathscr{D}(S/QS)$.
Indeed, when $\mathscr{A}$ is a $2$-arrangement,
it is shown in \cite{Nakashima-noeth} that
$\mathscr{D}(S/QS)$ is right and left Noetherian,
using free bases for $D^{(m)}(\mathscr{A})$ for all $m\geq 1$
constructed by Holm in \cite{Holm-phd}.
There are further results about $m$-freeness.
Coxeter arrangements of type A, B and D are $2$-free
(shown in \cite{N-Cox2free}).
We say that $\mathscr{A}$ is {\bf generic}
if $|\mathscr{A}|>\ell\geq 3$,
and if every $\ell$ hyperplanes of $\mathscr{A}$
intersect only at the origin.
For a generic arrangement $\mathscr{A}$, it is shown in \cite{NOS}
that $\mathscr{A}$ is $m$-free if and only if $m\geq|\mathscr{A}|-\ell+1$.
On the other hand, the behavior of $m$-freeness has
not been well analyzed yet when $m\geq 2$. Some basic questions remain open.
In particular, Holm in \cite{Holm-phd} asked the following two questions:
\begin{center}
\begin{itemize}
\item[Q1.] Does $m$-free imply $(m+1)$-free for any arrangement?
\item[Q2.] Are all arrangements $m$-free for $m$ large enough?
\end{itemize}
\end{center}
Q1 and Q2 are true for generic arrangements.
The aim of this paper is to prove that Q1 and Q2 are not true
by giving counter examples.
Our main results are as follows.
Proofs will appear in Section \ref{sec-proofs}.
\begin{Them}\label{thm-prod-arr}
Let $(\mathscr{A}_1,V_1)$ and $(\mathscr{A}_2,V_2)$ be arrangements
with $\dim(V_1)>0$ and $\dim(V_2)>0$.
The product arrangement
$(\mathscr{A}_1\times\mathscr{A}_2,V_1\oplus V_2)$ is $m$-free
if and only if both $(\mathscr{A}_1,V_1)$ and $(\mathscr{A}_2,V_2)$
are $i$-free for all $i\in\mathbb{N}$ with $1\leq i\leq m$.
Moreover, if $\exp_i(\mathscr{A}_1)=\{d_j^{(i)}\mid 1\leq j\leq s_i(\ell_1)\}$
and $\exp_{m-i}(\mathscr{A}_2)=\{e_k^{(m-i)}\mid 1\leq k\leq s_{m-i}(\ell_2)\}$
for $0\leq i\leq m$, then
$\exp_m(\mathscr{A}_1\times\mathscr{A}_2)=\bigcup_{i=0}^{m}
\{d_j^{(i)}+e_k^{(m-i)}\mid
1\leq j\leq s_i(\ell_1),\ 1\leq k\leq s_{m-i}(\ell_2)\}$.\hfill$\Box$
\end{Them}
Theorem \ref{thm-prod-arr} is the same as
\cite[Proposition 4.28]{Orlik-Terao} when $m=1$.
\begin{Cor}\label{cor-add_del14}
Let $(\mathscr{A}_{1},V_{1})$ and $(\mathscr{A}_{2},V_{2})$
be arrangements with $\dim(V_1)>0$ and $\dim(V_2)>0$.
The following are equivalent:
\begin{enumerate}
\item[(1)] $(\mathscr{A}_{1}\times\mathscr{A}_{2},V_{1}\oplus V_{2})$ is $m$-free.
\item[(2)] Both $(\mathscr{A}_{1},V_{1})$ and $(\mathscr{A}_{2},V_{2})$
are $i$-free for all $1\leq i\leq m$.
\item[(3)] $(\mathscr{A}_{1}\times\mathscr{A}_{2},V_{1}\oplus V_{2})$
is $i$-free for all $1\leq i\leq m$.\hfill$\Box$
\end{enumerate}
\end{Cor}
Theorem \ref{thm-prod-arr} (or Corollary \ref{cor-add_del14})
implies a counter example for Q2.
\begin{Examp}\label{nega-ans-Q1-redu}
{\rm
Let $\mathscr{A}$ and $\mathscr{A}^{\prime}$ be arrangements.
If $\mathscr{A}$ is not $1$-free, then the product arrangement
$\mathscr{A}\times\mathscr{A}^{\prime}$ is not $m$-free for any $m\geq 1$.
In particular, a generic arrangement is known
to be not $1$-free (see \cite{Orlik-Terao,Yuz}).
Hence if $\mathscr{A}$ is generic and if $\mathscr{A}^{\prime}$ is arbitrary,
then $\mathscr{A}\times\mathscr{A}^{\prime}$ is not $m$-free
for any $m\geq 1$.\hfill$\Box$
}
\end{Examp}
We say that $\mathscr{A}$ is {\bf reducible}
if $(\mathscr{A},V)=(\mathscr{A}_1,V_1)\times(\mathscr{A}_2,V_2)$
with $\dim(V_1)>0$ and $\dim(V_2)>0$ after a change of coordinates.
Otherwise $\mathscr{A}$ is said to be {\bf irreducible}.
Example \ref{nega-ans-Q1-redu} is of a reducible arrangement.
There also exists a counter example of an irreducible arrangement.
\begin{Prop}\label{nega-ans-Q1-irr}
Let $\mathscr{A}$ be a $4$-arrangement
defined by $Q=xyzw(x+y+z)(x+y+z+w)$.
Then $\mathscr{A}$ is not $m$-free for any $m\geq 1$.\hfill$\Box$
\end{Prop}
To prove Proposition \ref{nega-ans-Q1-irr},
we need some definitions and Proposition \ref{prop-A_X-free}.
Let
\begin{align}
L(\mathscr{A})=\left\{\bigcap_{H\in\mathscr{B}}H
\,\middle|\,\mathscr{B}\subseteq\mathscr{A}\right\}
\end{align}
be the set of all intersections of hyperplanes in $\mathscr{A}$,
which is partially ordered by the reverse inclusion.
We call $L(\mathscr{A})$ the {\bf intersection lattice} of $\mathscr{A}$.
\begin{Defi}
For $X\in L(\mathscr{A})$, a {\bf localization} $\mathscr{A}_X$
is a subarrangement of $\mathscr{A}$ defined by
\begin{align}
\mathscr{A}_X=\left\{H\in\mathscr{A}\,\middle|\,X\subseteq H\right\}.
\end{align}
\end{Defi}
\begin{Prop}\label{prop-A_X-free}
If $\mathscr{A}$ is $m$-free, then $\mathscr{A}_X$ is $m$-free
for all $X\in L(\mathscr{A})$.\hfill$\Box$
\end{Prop}
Proposition \ref{prop-A_X-free} is the same as
\cite[Theorem 4.37]{Orlik-Terao} when $m=1$.
The contraposition of Proposition \ref{prop-A_X-free} is useful to know that
arrangements are not $m$-free.

Next, we answer to Q1.
For $\ell\geq 2$, ${\rm Shi}_{\ell}$ is an $(\ell+1)$-arrangement defined by
\begin{align}\label{eq-shi-arr-A}
Q({\rm Shi}_{\ell})=
z\prod_{i=1}^{\ell}x_i(x_i-z)
\prod_{1\leq i<j\leq\ell}(x_i-x_j)(x_i-x_j-z).
\end{align}
Let $\Phi_{\ell}$ be the empty $\ell$-arrangement.
The arrangement ${\rm Shi}_{\ell}\times\Phi_1$
is the coning of a Shi arrangement (defined in \cite{Shi})
of the root system of the type A,
and ${\rm Shi}_{\ell}$ is known to be $1$-free
(see \cite{Athanasiadis, Yoshinaga-ER}).
For higher case, we have the following.
\begin{Them}\label{thm-shi-not2free}
The arrangement ${\rm Shi}_{\ell}$ is not $2$-free
for $\ell\geq 2$.\hfill$\Box$
\end{Them}
Therefore, ${\rm Shi}_{\ell}$ is an example such that
$1$-free does not imply $2$-free.
This means that Q1 is not true.

\section{Basic properties}
In this section, while we assume $m\geq 1$,
we introduce basic properties which are useful to observe $m$-freeness.
Although proofs of results in this section are already known,
we give their proofs, using the notations of this paper.
We first remark that the following relations hold:
\begin{align}\label{eq-relations-Weyl-alg}
x_i x_j=x_j x_i,\ \partial_i \partial_j=\partial_j \partial_i,\ 
\partial_i x_j=x_j \partial_i\ (i\neq j)\ {\rm and}\ 
\partial_i x_i=x_i \partial_i+1.
\end{align}
Indeed, the last relation follows from
$\partial_i x_i(f)=\partial_i(x_i f)=f+x_i \partial_i(f)=(1+x_i \partial_i)(f)$
$(f\in S)$, and the others are obvious.

\subsection{A criterion to know ideal stabilities}
For two operators $\theta,\eta\in {\rm End}_{\mathbb{K}}(S)$,
a commutator $[\theta,\eta]$ of $\theta$ and $\eta$ is defined by
$[\theta,\eta]=\theta\eta-\eta\theta$.
By the equation \eqref{eq-relations-Weyl-alg},
for $a\in\mathbb{N}\setminus\{0\}$,
\begin{align*}
[\partial_i^{a},x_i]=\partial_i^a x_i -x_i \partial_i^a
=x_i \partial_i^a +a\partial_i^{a-1} -x_i \partial_i^a=a\partial_i^{a-1}.
\end{align*}
Let $\bm{e}_i\in\mathbb{N}^{\ell}$ be the $i$th unit vector.
Then for $\bm{a}=(a_1\dots,a_{\ell})\in\mathbb{N}$
and for $1\leq i\leq\ell$,
\begin{align*}
[\partial^{\bm{a}},x_i]=\left\{
\begin{array}{ll}
a_i\partial^{\bm{a}-\bm{e}_i}&(a_i\neq 0),\\
0&(a_i=0).
\end{array}
\right.
\end{align*}
This implies that for $\theta\in D^{(m)}(S)$ and for $\alpha\in V^{\ast}$,
$[\theta,\alpha]\in D^{(m-1)}(S)$.
Moreover the following holds.
\begin{Lem}\label{lem-[]-in-D(A)}
For $\theta\in D^{(m)}(\mathscr{A})$ and for $\alpha\in V^{\ast}$,
the commutator $[\theta,\alpha]$ lies in $D^{(m-1)}(\mathscr{A})$.
\end{Lem}
\noindent
{\it Proof.}
Let $f\in S$. Then
$[\theta,\alpha](Qf)=\theta(\alpha Qf)-\alpha\theta(Qf)\in QS.$\hfill$\Box$
\vspace{4mm}

For an ideal $J$ of $S$, an $S$-submodule $D^{(m)}(J)$ of $D^{(m)}(S)$
is defined by
\begin{align*}
D^{(m)}(J)=\left\{\theta\in D^{(m)}(S)\,\middle|\,\theta(J)\subseteq J\right\}.
\end{align*}
\begin{Prop}[cf. Theorem 2.4 in \cite{Holm}]\label{prop-holm2.4}
\begin{align}
D^{(m)}(\mathscr{A})=\bigcap_{H\in\mathscr{A}}D^{(m)}(\alpha_H S).
\end{align}
\end{Prop}
\noindent
{\it Proof.}
Let $\theta\in \bigcap_{H\in\mathscr{A}}D^{(m)}(\alpha_H S)$.
For $H\in\mathscr{A}$,
$\theta(QS)\subseteq\theta(\alpha_H S)\subseteq \alpha_H S$.
Since the linear forms $\alpha_H$ $(H\in\mathscr{A})$ are coprime, we have
\begin{align*}
\theta(QS)\subseteq\left(\prod_{H\in\mathscr{A}}\alpha_H\right)S=QS.
\end{align*}
Thus the inclusion
$\bigcap_{H\in\mathscr{A}}D^{(m)}(\alpha_H S)\subseteq D^{(m)}(\mathscr{A})$
holds.

Next, we verify the converse inclusion
$D^{(m)}(\mathscr{A})\subseteq\bigcap_{H\in\mathscr{A}}D^{(m)}(\alpha_H S)$
by the double induction on $m$ and $|\mathscr{A}|$.
If $m=1$ then the assertion follows from the derivation case
(see \cite[Proposition 4.8]{Orlik-Terao}).
If $|\mathscr{A}|=1$ then it is obvious since $Q=\alpha_H$.

Suppose $m>1$ and $|\mathscr{A}|>1$. Let $\theta\in D^{(m)}(\mathscr{A})$.
Let $H_1$ be any hyperplane in $\mathscr{A}$, and let $\alpha_1=\alpha_{H_1}$.
By Lemma \ref{lem-[]-in-D(A)} and the induction hypothesis of $m$,
$[\theta,\alpha_1]\in D^{(m-1)}(\mathscr{A})
\subseteq \bigcap_{H\in\mathscr{A}}D^{(m-1)}(\alpha_H S)$.
Let $Q^{\prime}=\prod_{H\in\mathscr{A}\setminus\{H_1\}}\alpha_H$.
Then
$[\theta,\alpha_1](Q^{\prime}S)\subseteq Q^{\prime}S$
and $\theta(QS)\subseteq QS\subseteq Q^{\prime}S$.
For $f\in S$,
\begin{align}\label{eq-a1theta-DA-H1}
\alpha_1\theta(Q^{\prime}f)=\theta(Qf)-[\theta,\alpha_1](Q^{\prime}f)
\in Q^{\prime}S.
\end{align}
Since $\alpha_1$ and $Q^{\prime}$ are coprime,
the equation \eqref{eq-a1theta-DA-H1} implies that
$\theta$ belongs to $D^{(m)}(\mathscr{A}\setminus\{H_1\})$.
By the induction hypothesis of $|\mathscr{A}|$,
\begin{align*}
\theta\in \bigcap_{H\in\mathscr{A}\setminus\{H_1\}}D^{(m)}(\alpha_H S).
\end{align*}
Let $H_2\in\mathscr{A}$ satisfying $H_2\neq H_1$.
The same argument implies that
\begin{align*}
\theta\in\left(\bigcap_{H\in\mathscr{A}\setminus\{H_1\}}D^{(m)}(\alpha_H S)\right)
\bigcap\left(\bigcap_{H\in\mathscr{A}\setminus\{H_2\}}D^{(m)}(\alpha_H S)\right)
=\bigcap_{H\in\mathscr{A}}D^{(m)}(\alpha_H S).\quad\Box
\end{align*}

\begin{Lem}[cf. the proof of Proposition 2.10 in \cite{Singh}]\label{lem-singh}
Let $m\geq 2$ and let $H\in\mathscr{A}$. Then
$\theta\in D^{(m)}(\alpha_H S)$ if and only if
$[\theta,x_i]\in D^{(m-1)}(\alpha_H S)$
for all $i\in\mathbb{N}$ with $1\leq i\leq \ell$.
\end{Lem}
\noindent
{\it Proof.}
The assertion of `only if' is obvious.
Conversely we take $\theta\in D^{(m)}(S)$ such that
$[\theta,x_i](\alpha_H S)\subseteq\alpha_H S$
for all $i\in\mathbb{N}$ with $1\leq i\leq m$.
We verify that for any $\bm{a}=(a_1,\dots,a_{\ell})\in\mathbb{N}^{\ell}$,
$\theta(\alpha_H x^{\bm{a}})\in\alpha_H S$ by induction on $|\bm{a}|$.
If $|\bm{a}|=0$ then $\theta(\alpha_H x^{\bm{a}})=\theta(\alpha_H)=0$.
Let $|\bm{a}|\geq 1$.
Then there exist at least one index $i$ such that $a_i\neq 0$.
By the induction hypothesis, we obtain
\begin{align*}
\theta(\alpha_H x^{\bm{a}})=
[\theta,x_i](\alpha_H x^{\bm{a}-\bm{e}_i})+
x_i \theta(\alpha_H x^{\bm{a}-\bm{e}_i})\in\alpha_H S.\quad\Box
\end{align*}

\begin{Prop}[cf. Proposition 2.3 in \cite{Holm}]\label{prop-holm2.3}
Let $\theta\in D^{(m)}(S)$ and let $H\in\mathscr{A}$. Then
$\theta\in D^{(m)}(\alpha_H S)$ if and only if
$\theta(\alpha_H x^{\bm{a}})\in\alpha_H S$ for all $\bm{a}\in\mathbb{N}^{\ell}$
with $|\bm{a}|=m-1$.
\end{Prop}
\noindent
{\it Proof.}
The assertion of `only if' is obvious.
We verify the converse by induction on $m$.
If $m=1$ then the assumption (i.e., $\theta(\alpha_H)\in\alpha_H S$)
implies that for $f\in S$,
\begin{align*}
\theta(\alpha_H f)=\theta(\alpha_H)f+\alpha_H \theta(f)\in\alpha_H S.
\end{align*}

Let $m\geq 2$, $\bm{a}\in\mathbb{N}^{\ell}$
with $|\bm{a}|=m-2$ and let $i\in\mathbb{N}$ with $1\leq i\leq m$.
By the assumption,
$\theta(\alpha_H x^{\bm{a}+\bm{e}_i})\in\alpha_H S$.
Since $\deg(\alpha_H x^{\bm{a}})=m-1$,
we have $\theta(\alpha_H x^{\bm{a}})=0$. Then
\begin{align*}
[\theta,x_i](\alpha_H x^{\bm{a}})=
\theta(\alpha_H x^{\bm{a}+\bm{e}_i})-x_i \theta(\alpha_H x^{\bm{a}})
=\theta(\alpha_H x^{\bm{a}+\bm{e}_i})\in\alpha_H S.
\end{align*}
Therefore, by the induction hypothesis,
$[\theta,x_i]\in D^{(m-1)}(\alpha_H S)$.
This implies by Lemma \ref{lem-singh} that
$\theta\in D^{(m)}(\alpha_H S)$.
\hfill$\Box$
\vspace{4mm}

We summarize Proposition \ref{prop-holm2.4} and Proposition \ref{prop-holm2.3}
as follows.
\begin{Cor}\label{cor-check-D(A)}
Let $\mathscr{A}$ be an arrangement. Then
\begin{align}\label{eq-check-in-D(A)}
D^{(m)}(\mathscr{A})=\bigcap_{H\in\mathscr{A}}
\left\{\theta\in D^{(m)}(S)\,\middle|\,
\begin{array}{l}
\theta(\alpha_{H}x^{\bm{a}})\in \alpha_{H}S\ 
{\rm for\ all}\ \bm{a}\in\mathbb{N}^{\ell}\\
{\rm with}\ |\bm{a}|=m-1
\end{array}
\right\}.\quad\Box
\end{align}
\end{Cor}
Corollary \ref{cor-check-D(A)} is useful to know whether an operator
$\theta\in D^{(m)}(S)$ belongs to $D^{(m)}(\mathscr{A})$.
\begin{Examp}\label{ex-Euler}
{\rm
Let $\mathscr{A}$ be any arrangement.
The Euler operator
\begin{align}
\theta_E=\sum_{|\bm{a}|=m}\frac{m!}{\bm{a}!}x^{\bm{a}}\partial^{\bm{a}}
\end{align}
lies in $D^{(m)}(\mathscr{A})$ by Corollary \ref{cor-check-D(A)}.
Indeed, for $H\in\mathscr{A}$ and for $\bm{b}\in\mathbb{N}^{\ell}$
with $|\bm{b}|=m-1$,
\begin{align*}
\theta_E(\alpha_H x^{\bm{b}})=m!\alpha_H x^{\bm{b}}\in\alpha_H S.
\qquad\qquad\Box
\end{align*}
}
\end{Examp}
\begin{Examp}\label{ex-theta1theta2}
{\rm
Let $\mathscr{A}$ be a $2$-arrangement defined by $Q(\mathscr{A})=xy(x+y)$.
We denote $\partial_x=\frac{\partial}{\partial x}$
and $\partial_y=\frac{\partial}{\partial y}$.
Operators $\theta_1=x(x+y)\partial_x^2$ and $\theta_2=y(x+y)\partial_y^2$
lie in $D^{(2)}(\mathscr{A})$.
Indeed, since
\begin{align*}
&\theta_1(x\cdot x)=2x(x+y)\in xS,\ 
\theta_1(x\cdot y)=0\in xS,\ 
\theta_1(y\cdot x)=0\in xS,\\
&\theta_1(y\cdot y)=0\in yS,\ 
\theta_1((x+y)\cdot x)=2x(x+y)\in yS,\ 
\theta_1((x+y)\cdot y)=0\in yS,
\end{align*}
Corollary \ref{cor-check-D(A)} implies
$\theta_1\in D^{(2)}(\mathscr{A})$.
Similarly we have $\theta_2\in D^{(2)}(\mathscr{A})$.
\hfill$\Box$
}
\end{Examp}

\subsection{Saito's criterion}
Let
\begin{align}
s_{m}(\ell)=s_m=s=\binom{\ell+m-1}{m},\ 
t_{m}(\ell)=t_m=t=\binom{\ell+m-2}{m-1}.
\end{align}
We note that $t_m(\ell)=s_{m-1}(\ell)$. Let
\begin{align}
\Omega_m(\ell)=\Omega_m=\Omega=
\left\{\bm{a}\in\mathbb{N}^{\ell}\,\middle|\,|\bm{a}|=m\right\}.
\end{align}
We fix an ordering of $\Omega$ by $\Omega=\{\bm{a}(1),\dots,\bm{a}(s)\}$.
For $\theta_{1},\dots,\theta_{s}\in D^{(m)}(\mathscr{A})$,
a coefficient matrix $M_{m}(\theta_{1},\dots,\theta_{s_m})$
is an $s\times s$ matrix defined by
$M_{m}(\theta_{1},\dots,\theta_{s_m})=M_{m}=
\left(\frac{\theta_i\left(x^{\bm{a}}\right)}{\bm{a}!}
\right)_{1\leq i\leq s,\bm{a}\in\Omega}$.
In another description,
\begin{align}
M_m(\theta_{1},\dots,\theta_{s})=
\left(
\begin{matrix}
\frac{\theta_1\left(x^{\bm{a}(1)}\right)}{\bm{a}(1)!}&\cdots&\frac{\theta_s\left(x^{\bm{a}(1)}\right)}{\bm{a}(1)!}\\
\vdots& &\vdots\\
\frac{\theta_1\left(x^{\bm{a}(s)}\right)}{\bm{a}(s)!}&\cdots&\frac{\theta_s\left(x^{\bm{a}(s)}\right)}{\bm{a}(s)!}
\end{matrix}
\right).
\end{align}
\begin{Examp}\label{ex-coef-matrix}
{\rm
Let $\ell=2$ and let $m=2$. Then $s=3$ and $t=2$.
The coefficient matrix of
$\theta_E=x^2\partial_x^2+y^2\partial_y^2+2xy\partial_x \partial_y$,
$\theta_1=x(x+y)\partial_x^2$ and $\theta_2=y(x+y)\partial_y^2$
is the following:
\begin{align*}
M_m(\theta_E,\theta_1,\theta_2)=\left(
\begin{matrix}
x^2&x(x+y)&0\\
y^2&0&y(x+y)\\
2xy&0&0
\end{matrix}
\right).
\end{align*}
The rows are coefficients of $\partial_x^2$, $\partial_y^2$,
$\partial_x \partial_y$ from the top, and
the columns correspond to $\theta_E$, $\theta_1$, $\theta_2$
from the left.\hfill$\Box$
}
\end{Examp}
\begin{Prop}[Proposition I\hspace{-0.5mm}I\hspace{-0.5mm}I.5.2 in \cite{Holm-phd}]\label{add_del3}
If $\theta_{1},\dots,\theta_{s}\in D^{(m)}(\mathscr{A})$, then
\begin{align*}
\det M_{m}(\theta_{1},\dots,\theta_{s})\in Q^{t}S.
\end{align*}
\end{Prop}
\noindent
{\it Proof.}
Let $H\in\mathscr{A}$. We may assume $\alpha_H=x_1$. Since
\begin{align*}
|\{\bm{a}=(a_1,\dots,a_{\ell})\in\Omega_m\mid a_1\geq 1\}|
=|\Omega_{m-1}|=s_{m-1}=t,
\end{align*}
there exist $t$ rows in $M_m(\theta_{1},\dots,\theta_{s})$
such that all entries are divided by $x_1$.
Thus $\det M_m(\theta_{1},\dots,\theta_{s})\in\alpha_H^{t}S$.
Since $H$ is arbitrary,
$\det M_m(\theta_{1},\dots,\theta_{s})\in Q^{t}S$.
\hfill$\Box$
\vspace{4mm}

The following is Saito's criterion for $D^{(m)}(\mathscr{A})$
which is first given by Saito in \cite{Saito-criterion}
for $D^{(1)}(\mathscr{A})$ and which is generalized by
Holm in \cite{Holm-phd,Holm} for $D^{(m)}(\mathscr{A})$.
The proof is similar to \cite[Theorem 4.19]{Orlik-Terao}.
\begin{Them}[Saito's criterion, Proposition I\hspace{-0.5mm}I\hspace{-0.5mm}I.5.8 in \cite{Holm-phd}]\label{thm-saito's-criterion}
Given homogeneous operators
$\theta_{1},\dots,\theta_{s}\in D^{(m)}(\mathscr{A})$,
the following are equivalent:
\begin{enumerate}
\item[(1)] $\det M_{m}(\theta_{1},\dots,\theta_{s})
=cQ^t$ for some $c\in \mathbb{K}\setminus\{0\}$.
\item[(2)] $\theta_{1},\dots,\theta_{s}$ form a basis
for $D^{(m)}(\mathscr{A})$ over $S$.
\end{enumerate}
\end{Them}
\noindent
{\it Proof.}
$\left[(1)\ \Rightarrow\ (2)\right]$
Since $\det M_{m}(\theta_{1},\dots,\theta_{s})\neq 0$,
$\theta_{1},\dots,\theta_{s}$ are $S$-independent.
It is enough to prove that $\theta_{1},\dots,\theta_{s}$
generate $D^{(m)}(\mathscr{A})$ over $S$.
We may assume that $\det M_{m}(\theta_{1},\dots,\theta_{s})=Q^t$.
Let $\eta\in D^{(m)}(\mathscr{A})$.
Since $\theta_i=\sum_{\bm{a}\in\Omega}
\frac{\theta_i(x^{\bm{a}})}{\bm{a}!}\partial^{\bm{a}}$
for $1\leq i\leq s$,
we have by Cramer's rule that
\begin{align*}
Q^t\partial^{\bm{a}}\in S\theta_1+\cdots+S\theta_s
\end{align*}
for $\bm{a}\in\Omega$.
Then there exist $f_1,\dots,f_s\in S$ such that
\begin{align*}
Q^t\eta=f_1\theta_1+\cdots+f_s\theta_s.
\end{align*}
By Proposition \ref{add_del3},
$\det M_{m}(\theta_{1},\dots,\theta_{i-1},\eta,\theta_{i+1},\dots,\theta_{s})
\in Q^t S$.
Thus
\begin{align*}
Q^{2t} S&\ni
Q^t\det M_{m}(\theta_{1},\dots,\theta_{i-1},
\eta,\theta_{i+1},\dots,\theta_{s})\\
&=\det M_{m}(\theta_{1},\dots,\theta_{i-1},
Q^t\eta,\theta_{i+1},\dots,\theta_{s})\\
&=\det M_{m}(\theta_{1},\dots,\theta_{i-1},
f_i\theta_i,\theta_{i+1},\dots,\theta_{s})\\
&=f_iQ^t.
\end{align*}
This implies that $f_i\in Q^t S$ for $1\leq i\leq s$.
Therefore,
$\eta=(f_1/Q^t)\theta_1+\cdots+(f_s/Q^t)\theta_s
\in S\theta_1+\cdots+S\theta_s$.

$\left[(2)\ \Rightarrow\ (1)\right]$
By Proposition \ref{add_del3} and
the linear independence for $\theta_1,\dots,\theta_s$ over $S$,
$\det M_m(\theta_1,\dots,\theta_s)=fQ^t$
for some $f\in S\setminus\{0\}$.
Let $H\in\mathscr{A}$. We may assume that $\alpha_H=x_1$.
Also we may assume that
$\{\bm{a}(1),\dots,\bm{a}(t)\}=\{\bm{a}\in\Omega\mid a_1\geq 1\}$.
We define
\begin{align*}
\eta_i=\left\{
\begin{array}{ll}
Q\partial^{\bm{a}(i)}&{\rm if}\quad 1\leq i\leq t,\\
(Q/x_1)\partial^{\bm{a}(i)}&{\rm if}\quad t+1\leq i\leq s.
\end{array}
\right.
\end{align*}
Since $\eta_i\in S\theta_1+\cdots+S\theta_s$,
there exists an $s\times s$ matrix $N$ whose entries lie in $S$ such that
$M_m(\eta_1,\dots.\eta_s)=M_m(\theta_1,\dots,\theta_s)N$.
Then
\begin{align*}
Q^s/x_1^{s-t}=\det M_m(\eta_1,\dots.\eta_s)
=\det M_m(\theta_1,\dots,\theta_s)\det N
=fQ^t \det N.
\end{align*}
This implies that $(Q/\alpha_H)^{s-t}=f\det N$, and therefore
$f$ divides $(Q/\alpha_H)^{s-t}$.
This is true for all $H\in\mathscr{A}$.
Since polynomials $\{(Q/\alpha_H)^{s-t}\}_{H\in\mathscr{A}}$
have no common factor,
we obtain $f\in\mathbb{K}\setminus\{0\}$.
\hfill$\Box$
\vspace{4mm}

We note that Theorem \ref{thm-saito's-criterion} (1) does not depend on
the ordering of $\Omega$. In other words,
to use Theorem \ref{thm-saito's-criterion}, we may choose any ordering
of $\Omega$.

\begin{Them}[Proposition I\hspace{-0.5mm}I\hspace{-0.5mm}I.5.9 in \cite{Holm-phd}]\label{add_del6}
Let $\theta_{1},\dots,\theta_{s}\in D^{(m)}(\mathscr{A})$
be homogeneous operators which are $S$-independent.
Then $\mathscr{A}$ is $m$-free with homogeneous basis
$\theta_{1},\dots,\theta_{s_{m}}$ if and only if
$\sum_{j=1}^{s}\deg\left(\theta_{j}\right)=t|\mathscr{A}|$.\hfill$\Box$
\end{Them}
\begin{Examp}\label{ex-elementary-saito-criterion}
{\rm
Let $\ell=2$, $m=2$ and let $Q(\mathscr{A})=xy(x+y)$.
Then $s=3$ and $t=2$.
Operators
$\theta_E=x^2\partial_x^2+y^2\partial_y^2+2xy\partial_x \partial_y$,
$\theta_1=x(x+y)\partial_x^2$ and $\theta_2=y(x+y)\partial_y^2$
lie in $D^{(2)}(\mathscr{A})$ by Example \ref{ex-Euler}
and Example \ref{ex-theta1theta2}. Since
\begin{align}
\det M_m(\theta_E,\theta_1,\theta_2)=\left|
\begin{matrix}
x^2&x(x+y)&0\\
y^2&0&y(x+y)\\
2xy&0&0
\end{matrix}
\right|=2Q^2,
\end{align}
the operators $\theta_E,\theta_1,\theta_2$
form a basis for $D^{(2)}(\mathscr{A})$.\hfill$\Box$
}
\end{Examp}
\begin{Examp}[Proposition I\hspace{-0.5mm}I\hspace{-0.5mm}I.6.7 in \cite{Holm-phd}]\label{ex-ell=2basis}
{\rm
Let $\ell=2$ and let $\mathscr{A}=\{H_1,\dots,H_n\}$, where $n=|\mathscr{A}|$.
We may assume that
\begin{align*}
H_1=\{x=0\},\quad H_j=\{a_j x+y=0\}\quad(2\leq j\leq n)
\end{align*}
for some distinct scalars $a_2,\dots,a_r\in\mathbb{K}$.
We define $Q_1=Q/x$ and $Q_j=Q/(a_i x+y)$
for $2\leq j\leq n$. Then
\begin{align*}
Q_1\partial_{y}^m,\quad Q_j(\partial_{x}-a_j \partial_{y})^m
\quad(2\leq j\leq n)
\end{align*}
belong to $D^{(m)}(\mathscr{A})$.
By Theorem \ref{add_del6}, $\mathscr{A}$ is $m$-free
with the following basis:
\begin{enumerate}
\item[(1)] $\{\theta_E,Q_1\partial_{y}^m,Q_2(\partial_{x}-a_2 \partial_{y})^m,
\dots,Q_m(\partial_{x}-a_m \partial_{y})^m\}$ if $m\leq n-2$.
\item[(2)] $\{Q_1\partial_{y}^m,Q_2(\partial_{x}-a_2 \partial_{y})^m,
\dots,Q_n(\partial_{x}-a_n \partial_{y})^m\}$ if $m=n-1$.
\item[(3)] $\{Q_1\partial_{y}^m,Q_2(\partial_{x}-a_2 \partial_{y})^m,
\dots,Q_n(\partial_{x}-a_n \partial_{y})^m,Q\eta_{n+1},\dots,Q\eta_{m+1}\}$
if $m\geq n$, where
$\{\partial_{y}^m,(\partial_{x}-a_2 \partial_{y})^m,\dots,
(\partial_{x}-a_n \partial_{y})^m,\eta_{n+1},\dots,\eta_{m+1}\}$
is a basis for\\
$\sum_{j=0}^{m}\mathbb{K}\partial_{x}^j\partial_{y}^{m-j}$
over $\mathbb{K}$.\hfill$\Box$
\end{enumerate}
}
\end{Examp}
Example \ref{ex-elementary-saito-criterion} is a special case of
Example \ref{ex-ell=2basis}.

\section{Proofs}\label{sec-proofs}
\subsection{Proofs of Theorem \ref{thm-prod-arr} and Proposition \ref{nega-ans-Q1-irr}}
Let $(\mathscr{A}_1,V_1)$ and $(\mathscr{A}_2,V_2)$ be arrangements
with $\ell_1=\dim(V_1)>0$ and $\ell_2=\dim(V_2)>0$.
Let $S_1=\Sym(V_1^{\ast})=\mathbb{K}[x_1,\dots,x_{\ell_1}]$ and
$S_2=\Sym(V_2^{\ast})=\mathbb{K}[y_1,\dots,y_{\ell_2}]$.
We denote $Q_1=Q(\mathscr{A}_1)\in S_1$ and $Q_2=Q(\mathscr{A}_2)\in S_2$.
Let $\mathscr{A}=\mathscr{A}_1\times\mathscr{A}_2$ and
$S=S_1\otimes S_2=\mathbb{K}[x_1,\dots,x_{\ell_1},y_1,\dots,y_{\ell_2}]$.
Then $Q=Q(\mathscr{A})=Q_1 Q_2\in S$.

For $\theta\in D^{(i)}(S_1)$ and for $\eta\in D^{(j)}(S_2)$,
the product $\theta\eta=\eta\theta$ is commutative
in ${\rm End}_{\mathbb{K}}(S)$, by the equation \eqref{eq-relations-Weyl-alg}.
Let $SD^{(i)}(\mathscr{A}_1)D^{(j)}(\mathscr{A}_2)$ be
the $S$-submodule of ${\rm End}_{\mathbb{K}}(S)$ generated by
$\{\theta\eta\mid \theta\in D^{(i)}(\mathscr{A}_1),\ 
\eta\in D^{(j)}(\mathscr{A}_2)\}$.
We note that $D^{(m)}(\Phi_{\ell})=D^{(m)}(S)$.
\begin{Lem}
$SD^{(i)}(\mathscr{A}_1)D^{(j)}(\mathscr{A}_2)
\subseteq D^{(i+j)}(\mathscr{A}_1\times \mathscr{A}_2)$.
\end{Lem}
\noindent
{\it Proof.}
For $\theta\in D^{(i)}(\mathscr{A}_1)$, $\eta\in D^{(j)}(\mathscr{A}_2)$
and $f\in S$,
\begin{align*}
\theta\eta(Qf)=\sum_{k=1}^{K}\sum_{l=1}^L\theta(Q_1 f_k)\eta(Q_2 g_l)
\in Q_1 S_1 Q_2 S_2=QS,
\end{align*}
where $f=\sum_{k=1}^K\sum_{l=1}^L f_k g_l$ for some $f_k\in S_1$, $g_l\in S_2$.
\hfill$\Box$
\vspace{4mm}

We recall that $D^{(0)}(S_1)=S_1$ and $D^{(0)}(S_2)=S_2$.
The $S$-module $D^{(m)}(S)$ is decomposed as follows:
\begin{align}\label{eq-docompese-D(S)}
D^{(m)}(S)=\bigoplus_{i=0}^{m}SD^{(i)}(S_{1})D^{(m-i)}(S_{2}).
\end{align}
In this section, a partial derivative of $x_i$
is denoted by $\partial_i=\frac{\partial}{\partial x_i}$
and that of $y_i$ is denoted by $\delta_i=\frac{\partial}{\partial y_i}$.
We note that $D^{(0)}(\mathscr{A}_1)=S_1$ and $D^{(0)}(\mathscr{A}_2)=S_2$.

\begin{Lem}\label{lem-for-decomp1}
Let $\mathscr{A}=\mathscr{A}_1\times\mathscr{A}_2$.
Then the $S$-module $D^{(m)}(\mathscr{A})$ is decomposed as follows:
\begin{align*}
D^{(m)}(\mathscr{A})=\bigoplus_{i=0}^{m}
\left(SD^{(i)}(\mathscr{A}_{1})D^{(m-i)}(S_{2})\cap
SD^{(i)}(S_{1})D^{(m-i)}(\mathscr{A}_{2})\right).
\end{align*}
\end{Lem}
\noindent
{\it Proof.}
Let $\theta\in SD^{(i)}(\mathscr{A}_{1})D^{(m-i)}(S_{2})\cap
SD^{(i)}(S_{1})D^{(m-i)}(\mathscr{A}_{2})$ and let $f\in S$.
Since $\theta\in SD^{(i)}(\mathscr{A}_{1})D^{(m-i)}(S_{2})$,
we have $\theta(Qf)\in Q_1 S$.
Similarly $\theta(Qf)\in Q_2 S$.
Since $Q_1$ and $Q_2$ are coprime, we have $\theta(Qf)\in QS$.
This means $\theta\in D^{(m)}(\mathscr{A})$.

Conversely let $\theta\in D^{(m)}(\mathscr{A})$.
By the equation \eqref{eq-docompese-D(S)}, we can describe that
$\theta=\sum_{i=0}^{m}\theta^{(i)}$ for some
$\theta^{(i)}\in SD^{(i)}(S_{1})D^{(m-i)}(S_{2})$.
By the symmetry of $\mathscr{A}_1$ and $\mathscr{A}_2$,
it is enough to verify that
$\theta^{(i)}\in D^{(i)}(\mathscr{A}_{1})D^{(m-i)}(S_{2})$ for $0\leq i\leq m$.
If $i=0$ then we immediately have
$\theta^{(0)}\in SD^{(0)}(S_{1})D^{(m)}(S_{2})
=SD^{(0)}(\mathscr{A}_{1})D^{(m)}(S_{2})$.
Let $i>0$. Since $\{y^{\bm{a}}\delta^{\bm{b}}\mid
\bm{a},\bm{b}\in\mathbb{N}^{\ell_2},|\bm{b}|=m-i\}$
is a $\mathbb{K}$-basis for $D^{(m-i)}(S_{2})$,
the operator $\theta^{(i)}$ is described as
\begin{align*}
\theta^{(i)}=\sum_{\bm{a},\bm{b}\in \mathbb{N}^{\ell_2},|\bm{b}|=m-i}
y^{\bm{a}}\delta^{\bm{b}}\theta_{\bm{a},\bm{b}}
\end{align*}
for some $\theta_{\bm{a},\bm{b}}\in D^{(i)}(S_{1})$.
Now let $H\in\mathscr{A}_{1}$,
$\bm{c}\in\mathbb{N}^{\ell_1}$ with $|\bm{c}|=i-1$ and let
$\bm{d}\in\mathbb{N}^{\ell_2}$ with $|\bm{d}|=m-i$.
Since $\theta^{(j)}(\alpha_{H}x^{\bm{c}}y^{\bm{d}})=0$ for $j\neq i$,
we have
\begin{align*}
\alpha_{H}S\ni \theta(\alpha_{H}x^{\bm{c}}y^{\bm{d}})
=\theta^{(i)}(\alpha_{H}x^{\bm{c}}y^{\bm{d}})
=\sum_{\bm{a}\in\mathbb{N}^{\ell_2}}y^{\bm{a}}\bm{d}!
\theta_{\bm{a},\bm{d}}(\alpha_{H}x^{\bm{c}}),
\end{align*}
i.e., there exists $f\in S$ such that
$\alpha_H f=\sum_{\bm{a}\in\mathbb{N}^{\ell_2}}y^{\bm{a}}\bm{d}!
\theta_{\bm{a},\bm{d}}(\alpha_{H}x^{\bm{c}})$.
Here $\{y^{\bm{a}}\mid\bm{a}\in\mathbb{N}^{\ell_2}\}$
is an $S_1$-basis for $S$.
Then there exist $f_{\bm{a}}\in S_1$ $(\bm{a}\in\mathbb{N}^{\ell_2})$ such that
\begin{align*}
\sum_{\bm{a}\in\mathbb{N}^{\ell_2}}y^{\bm{a}}
\alpha_H f_{\bm{a}}
=\sum_{\bm{a}\in\mathbb{N}^{\ell_2}}y^{\bm{a}}\bm{d}!
\theta_{\bm{a},\bm{d}}(\alpha_{H}x^{\bm{c}}).
\end{align*}
Since $\theta_{\bm{a},\bm{d}}(\alpha_{H}x^{\bm{c}})\in S_1$, we have
$\theta_{\bm{a},\bm{d}}(\alpha_{H}x^{\bm{c}})
=\alpha_H f_{\bm{a}}\in\alpha_H S_1$
for any $\bm{a}\in\mathbb{N}^{\ell_2}$.
This is true for any $\bm{c}\in\mathbb{N}^{\ell_1}$ with $|\bm{c}|=i-1$
and for any $\bm{d}\in\mathbb{N}^{\ell_2}$ with $|\bm{d}|=m-i$.
By Corollary \ref{cor-check-D(A)},
$\theta_{\bm{a},\bm{d}}\in D^{(i)}(\mathscr{A}_1)$
for any $\bm{a},\bm{d}\in\mathbb{N}^{\ell_2}$ with $|\bm{d}|=m-i$.
Hence we conclude
$\theta^{(i)}\in SD^{(i)}(\mathscr{A}_{1})D^{(m-i)}(S_{2})$.
\hfill$\Box$

\begin{Lem}\label{lem-for-decomp2}
For $0\leq i\leq m$,
\begin{align*}
SD^{(i)}(\mathscr{A}_{1})D^{(m-i)}(\mathscr{A}_{2})=
SD^{(i)}(\mathscr{A}_{1})D^{(m-i)}(S_{2})
\cap SD^{(i)}(S_{1})D^{(m-i)}(\mathscr{A}_{2}).
\end{align*}
\end{Lem}
\noindent
{\it Proof.}
It is obvious that $SD^{(i)}(\mathscr{A}_{1})D^{(m-i)}(\mathscr{A}_{2})
\subseteq SD^{(i)}(\mathscr{A}_{1})D^{(m-i)}(S_{2})
\cap SD^{(i)}(S_{1})$\\
$D^{(m-i)}(\mathscr{A}_{2})$.

Conversely let
$\theta\in SD^{(i)}(\mathscr{A}_{1})D^{(m-i)}(S_{2})
\cap SD^{(i)}(S_{1})D^{(m-i)}(\mathscr{A}_{2})$.
Let $\{\eta_{\lambda}\}_{\lambda\in\Lambda}$
be a $\mathbb{K}$-basis for $D^{(m-i)}(\mathscr{A}_{2})$.
Since $SD^{(i)}(S_{1})D^{(m-i)}(\mathscr{A}_{2})
=D^{(i)}(S_{1})\otimes_{\mathbb{K}}D^{(m-i)}(\mathscr{A}_{2})$,
there exist
$\theta_{\lambda}\in D^{(i)}(S_1)$ $(\lambda\in\Lambda)$ such that
$\theta=\sum_{\lambda\in\Lambda}\theta_{\lambda}\eta_{\lambda}$.
We show that
$\theta_{\lambda}\in D^{(i)}(\mathscr{A}_1)$ for any $\lambda\in\Lambda$,
i.e., $\theta_{\lambda}(Q_1 f)\in Q_1 S_1$ for any $f\in S_1$.
An operator $\theta_{0}$ is defined by
\begin{align*}
\theta_{0}=\sum_{\lambda\in\Lambda}\theta_{\lambda}(Q_1 f)\eta_{\lambda}
\in SD^{(m-i)}(\mathscr{A}_2)
\end{align*}
for any $f\in S_1$.
We note that $\theta_{0}(g)=\theta(Q_1 fg)$ for any $g\in S_2$.
Since $\{\delta^{\bm{b}}\mid\bm{b}\in\mathbb{N}^{\ell_2},|\bm{b}|=m-i\}$
is an $S$-basis for $SD^{(m-i)}(S_2)$, we can describe that
\begin{align*}
\theta_{0}=\sum_{\bm{b}\in\mathbb{N}^{\ell_2},|\bm{b}|=m-i}
f_{\bm{b}}\delta^{\bm{b}}
\end{align*}
for some $f_{\bm{b}}\in S$ $(\bm{b}\in\mathbb{N}^{\ell_2})$.
By the assumption $\theta\in D^{(i)}(\mathscr{A}_{1})D^{(m-i)}(S_{2})$,
for any $\bm{b}\in\mathbb{N}^{\ell_2}$ with $|\bm{b}|=m-i$,
\begin{align*}
f_{\bm{b}}=\frac{\theta_{0}(y^{\bm{b}})}{\bm{b}!}
=\frac{\theta(Q_1 f y^{\bm{b}})}{\bm{b}!}\in Q_1 S.
\end{align*}
This implies
\begin{align*}
\sum_{\lambda\in\Lambda}\theta_{\lambda}(Q_1 x^{\bm{a}})\eta_{\lambda}
=\theta_{0}\in Q_1 SD^{(m-i)}(S_2)\cap SD^{(m-i)}(\mathscr{A}_2).
\end{align*}
Here we prove $Q_1 SD^{(m-i)}(S_2)\cap SD^{(m-i)}(\mathscr{A}_2)
=Q_1 SD^{(m-i)}(\mathscr{A}_2)$.
Indeed, we assume that $\theta\in SD^{(m-i)}(S_2)$ and
$Q_1 \theta\in SD^{(m-i)}(\mathscr{A}_2)$.
Then we can describe that
$\theta=\sum_{\bm{a}\in\mathbb{N}^{\ell_1}}x^{\bm{a}}\theta_{\bm{a}}$
for some $\theta_{\bm{a}}\in D^{(m-i)}(S_2)$
$(\bm{a}\in\mathbb{N}^{\ell_1})$, while $Q_1 \theta(Q_2g)\in Q_2 S$
for any $g\in S_2$.
Since $Q_1$ and $Q_2$ are coprime, we have $\theta(Q_2g)\in Q_2 S$.
Thus there exist $g_{\bm{a}}\in S_2$ $(\bm{a}\in\mathbb{N}^{\ell_1})$ such that
$\sum_{\bm{a}\in\mathbb{N}^{\ell_1}}x^{\bm{a}}\theta_{\bm{a}}(Q_2 g)
=\sum_{\bm{a}\in\mathbb{N}^{\ell_1}}x^{\bm{a}} Q_2 g_{\bm{a}}$, that is,
$\theta_{\bm{a}}(Q_2 g)\in Q_2 S_2$ for any $\bm{a}\in\mathbb{N}^{\ell_1}$.
This means $\theta_{\bm{a}}\in D^{(m-i)}(\mathscr{A}_2)$, and hence
$Q_1 \theta\in Q_1 SD^{(m-i)}(\mathscr{A}_2).$
As a result, we have
\begin{align*}
\sum_{\lambda\in\Lambda}\theta_{\lambda}(Q_1 x^{\bm{a}})\eta_{\lambda}
\in Q_1 SD^{(m-i)}(\mathscr{A}_2).
\end{align*}

Since $S_1$ is flat over $\mathbb{K}$ and
$\eta_{\lambda}$ $(\lambda\in\Lambda)$ are $S_1$-torsion free elements,
$\{\eta_{\lambda}\}_{\lambda\in\Lambda}$ is an $S_1$-basis for
$SD^{(m-i)}(\mathscr{A}_{2})=S_1\otimes_{\mathbb{K}}
D^{(m-i)}(\mathscr{A}_{2})$.
Then there exist $f_{\lambda}\in S_1$ $(\lambda\in\Lambda)$ such that
\begin{align*}
\sum_{\lambda\in\Lambda}\theta_{\lambda}(Q_1 x^{\bm{a}})\eta_{\lambda}
=\sum_{\lambda\in\Lambda}Q_1 f_{\lambda}\eta_{\lambda},
\end{align*}
while $\theta_{\lambda}(Q_1 f)=Q_1 f_{\lambda}\in Q_1 S_1$
for any $\lambda\in\Lambda$.
\hfill$\Box$
\vspace{4mm}

Lemma \ref{lem-for-decomp1} and Lemma \ref{lem-for-decomp2}
imply the following.
\begin{Prop}\label{prop-decomp-D(A1timesA2)}
Let $\mathscr{A}=\mathscr{A}_1\times\mathscr{A}_2$. Then
\begin{align}
D^{(m)}(\mathscr{A})=\bigoplus_{i=0}^{m}
SD^{(i)}(\mathscr{A}_{1})D^{(m-i)}(\mathscr{A}_{2}).\qquad\Box
\end{align}
\end{Prop}
\begin{Examp}
{\rm
We consider an example of $\ell_1=2$, $\ell_2=1$ and $m=2$.
Let $\mathscr{A}_1$ be a $2$-arrangement defined by $Q(\mathscr{A}_1)=xy(x+y)$
and $\mathscr{A}_2$ the empty $1$-arrangement with the coordinate $z$.
Let $\mathscr{A}=\mathscr{A}_1\times\mathscr{A}_2$.
The sets
$\{\theta_E^{(2)}=x^2\partial_x^2+y^2\partial_y^2+2xy\partial_x \partial_y,\ 
\theta_1^{(2)}=x(x+y)\partial_x^2,\ \theta_2^{(2)}=y(x+y)\partial_y^2\}$
and
$\{\theta_E^{(1)}=x\partial_x+y\partial_y,\ \theta_1^{(1)}=y(x+y)\partial_y\}$
are bases for $D^{(2)}(\mathscr{A}_1)$ and $D^{(1)}(\mathscr{A}_1)$,
respectively.
Then the exponents are $\exp_{2}(\mathscr{A}_1)=\{2,\ 2,\ 2\}$ and
$\exp_{1}(\mathscr{A}_1)=\{1,\ 2\}$.
Meanwhile the sets $\{\partial_z^2\}$ and $\{\partial_z\}$
are bases for $D^{(2)}(\mathscr{A}_2)$ and $D^{(1)}(\mathscr{A}_2)$,
respectively.
The exponents are $\exp_{2}(\mathscr{A}_2)=\{0\}$ and
$\exp_{1}(\mathscr{A}_2)=\{0\}$. Then
Proposition \ref{prop-decomp-D(A1timesA2)} implies that
\begin{align*}
D^{(2)}(\mathscr{A})&=
SD^{(0)}(\mathscr{A}_{1})D^{(2)}(\mathscr{A}_{2})\oplus
SD^{(1)}(\mathscr{A}_{1})D^{(1)}(\mathscr{A}_{2})\oplus
SD^{(2)}(\mathscr{A}_{1})D^{(0)}(\mathscr{A}_{2})\\
&=S\partial_z^2\oplus S\theta_E^{(1)}\partial_z\oplus S\theta_1^{(1)}\partial_z
\oplus S\theta_E^{(2)}\oplus S\theta_1^{(2)}\oplus S\theta_2^{(2)}.
\end{align*}
Therefore, $\exp_{2}(\mathscr{A})=\{0,\ 1,\ 2,\ 2,\ 2,\ 2\}$.\hfill$\Box$
}
\end{Examp}

By Proposition \ref{prop-decomp-D(A1timesA2)},
we can prove Theorem \ref{thm-prod-arr}.
\vspace{4mm}

\noindent
{\it Proof of Theorem \ref{thm-prod-arr}.}
Let $(\mathscr{A}_{1},V_{1})$ and $(\mathscr{A}_{2},V_{2})$
be $i$-free for all $1\leq i\leq m$.
Let $\{\theta_{j}^{(i)}\mid 1\leq j\leq s_{i}(\ell_1)\}$
and $\{\eta_{k}^{(m-i)}\mid 1\leq k\leq s_{m-i}(\ell_{2})\}$
be homogeneous bases for $D^{(i)}(\mathscr{A}_{1})$ and
$D^{(m-i)}(\mathscr{A}_{2})$ for all $0\leq i\leq m$, respectively.
We note that $\{1\}$ is a basis for $D^{(0)}(\mathscr{A}_{1})$ over $S_1$
and for $D^{(0)}(\mathscr{A}_{2})$ over $S_2$.

We set $B_i=\{\theta_{j}^{(i)}\eta_{k}^{(m-i)}\mid
1\leq j\leq s_{i}(\ell_1),\ 1\leq k\leq s_{m-i}(\ell_{2})\}$.
To verify that $\bigcup_{i=0}^{m}B_i$ is linearly independent over $S$,
it suffices to show that
$B_i$ is linearly independent over $S$ for $0\leq i\leq m$.
We assume that
\begin{align}\label{eq-lin-indep-DA1DA2}
\sum_{j=1}^{s_{i}(\ell_1)}\,\sum_{k=1}^{s_{m-i}(\ell_2)}
f_{j,k}\theta_{j}^{(i)}\eta_{k}^{(m-i)}=0,\quad
{\rm where}\ f_{j,k}\in S.
\end{align}
Let $\bm{a}\in \mathbb{N}^{\ell_1}$ with $|\bm{a}|=i$.
By substituting $x^{\bm{a}}$ for all $\theta_{j}^{(i)}$
in the equation \eqref{eq-lin-indep-DA1DA2},
we have the following equation of operators
in $SD^{(m-i)}(S_2)$:
\begin{align*}
\sum_{k=1}^{s_{m-i}(\ell_2)}\left(\sum_{j=1}^{s_{i}(\ell_1)}
f_{j,k}\theta_{j}^{(i)}(x^{\bm{a}})\right)\eta_{k}^{(m-i)}=0.
\end{align*}
Since $S$ is flat over $S_2$ and
$\eta_{k}^{(m-i)}$ $(1\leq k\leq s_{m-i}(\ell_{2}))$
are $S$-torsion free elements,
$\{\eta_{k}^{(m-i)}\mid 1\leq k\leq s_{m-i}(\ell_{2})\}$
is linearly independent over $S$.
Then $\sum_{j=1}^{s_{i}(\ell_1)}f_{j,k}\theta_{j}^{(i)}(x^{\bm{a}})=0$
for $1\leq k\leq s_{m-i}(\ell_{2})$ and 
for $\bm{a}\in \mathbb{N}^{\ell_1}$ with $|\bm{a}|=i$.
This implies $\sum_{j=1}^{s_{i}(\ell_1)}f_{j,k}\theta_{j}^{(i)}=0$
for $1\leq k\leq s_{m-i}(\ell_{2})$.
Since $\{\theta_{j}^{(i)}\mid 1\leq j\leq s_{i}(\ell_1)\}$
is linearly independent over $S$, we have $f_{j,k}=0$
for $1\leq j\leq s_{i}(\ell_1)$ and for $1\leq k\leq s_{m-i}(\ell_{2})$.
Hence we conclude that $\bigcup_{i=0}^{m}B_i$
is linearly independent over $S$.

Next, we verify that the sum of degrees of $\bigcup_{i=0}^{m}B_i$
is equal to $t_m(\ell_1+\ell_2)|\mathscr{A}_1\times\mathscr{A}_2|$ as follows:
\begin{align*}
&\sum_{i=0}^{m}\sum_{j=1}^{s_{i}(\ell_1)}\,\sum_{k=1}^{s_{m-i}(\ell_2)}
\deg\left(\theta_{j}^{(i)}\eta_{k}^{(m-i)}\right)\\
=&\,\sum_{i=0}^{m}\sum_{j=1}^{s_{i}(\ell_1)}\,\sum_{k=1}^{s_{m-i}(\ell_2)}
\left(\deg\theta_{j}^{(i)}+\deg\eta_{k}^{(m-i)}\right)\\
=&\,\sum_{i=0}^{m}\left(
s_{m-i}(\ell_2)\sum_{j=1}^{s_{i}(\ell_1)}\deg\theta_{j}^{(i)}+
s_{i}(\ell_1)\sum_{k=1}^{s_{m-i}(\ell_2)}\deg\eta_{k}^{(m-i)}\right)\\
=&\,\sum_{i=1}^{m}s_{m-i}(\ell_2)t_i(\ell_1)|\mathscr{A}_1|+
\sum_{i=0}^{m-1}s_{i}(\ell_1)t_{m-i}(\ell_2)|\mathscr{A}_2|\\
=&\,\sum_{i=0}^{m-1}s_{m-1-i}(\ell_2)s_i(\ell_1)|\mathscr{A}_1|+
\sum_{i=0}^{m-1}s_{i}(\ell_1)s_{m-1-i}(\ell_2)|\mathscr{A}_2|\\
=&\,s_{m-1}(\ell_1+\ell_2)\left(|\mathscr{A}_1|+|\mathscr{A}_2|\right)\\
=&\,t_m(\ell_1+\ell_2)|\mathscr{A}_1\times\mathscr{A}_2|.
\end{align*}
Therefore, $\bigcup_{i=0}^{m}B_i$ is a basis for
$D^{(m)}(\mathscr{A}_1\times\mathscr{A}_2)$ over $S$ by Theorem \ref{add_del6}.

Conversely let $(\mathscr{A}_1\times\mathscr{A}_2,V_1\oplus V_2)$
be $m$-free. Let $\{\theta_{\lambda}^{(i)}\}_{\lambda\in\Lambda}$
and $\{\eta_{\mu}^{(m-i)}\}_{\mu\in M}$
be minimal sets of homogeneous generators for
$D^{(i)}(\mathscr{A}_1)$ over $S_1$ and
for $D^{(m-i)}(\mathscr{A}_2)$ over $S_2$, respectively.
Obviously,
$\{\theta_{\lambda}^{(i)}\eta_{\mu}^{(m-i)}\}_{\lambda\in\Lambda,\mu\in M}$
generates $D^{(i)}(\mathscr{A}_1)D^{(m-i)}(\mathscr{A}_2)$ over $S$.
We assume that
$B_i=\{\theta_{\lambda}^{(i)}\eta_{\mu}^{(m-i)}\}_{\lambda\in\Lambda,\mu\in M}$
is not minimal. Then there exist indices $\lambda_0,\ \mu_0$
and polynomials $f_{\lambda,\mu}\in S\ (\lambda\neq\lambda_0, \mu\neq\mu_0)$
such that
\begin{align}\label{eq-rela-notminimal}
\theta_{\lambda_0}^{(i)}\eta_{\mu_0}^{(m-i)}=
\sum_{\lambda\neq\lambda_0,\ \mu\neq\mu_0}
f_{\lambda,\mu}\theta_{\lambda}^{(i)}\eta_{\mu}^{(m-i)}.
\end{align}
If $\theta_{\lambda_0}^{(i)}\left(x^{\bm{a}}\right)=0$
for all $\bm{a}\in\mathbb{N}^{\ell_{1}}$
with $|\bm{a}|=i$, then $\theta_{\lambda_0}^{(i)}=0$.
This contradicts to the minimality of
$\{\theta_{\lambda}^{(i)}\}_{\lambda\in\Lambda}$.
Thus there exists $\bm{a}\in\mathbb{N}^{\ell_{1}}$ with $|\bm{a}|=i$
such that $\theta_{\lambda_0}^{(i)}\left(x^{\bm{a}}\right)\neq 0$.
By substituting $x^{\bm{a}}$ for all
$\theta_{\lambda}^{(i)}$ $(\lambda\in\Lambda)$
in the equation \eqref{eq-rela-notminimal},
\begin{align*}
\theta_{\lambda_0}^{(i)}\left(x^{\bm{a}}\right)\eta_{\mu_0}^{(m-i)}=
\sum_{\lambda\neq\lambda_0,\ \mu\neq\mu_0}
f_{\lambda,\mu}\theta_{\lambda}^{(i)}\left(x^{\bm{a}}\right)\eta_{\mu}^{(m-i)}.
\end{align*}
Since $\theta_{\lambda_0}^{(i)}\left(x^{\bm{a}}\right)\in S_1$ and
$f_{\lambda,\mu}\theta_{\lambda}^{(i)}\left(x^{\bm{a}}\right)\in S$,
we can describe that
\begin{align*}
\sum_{\bm{c}\in\mathbb{N}^{\ell_1}}x^{\bm{c}}g_{0,\bm{c}}
\eta_{\mu_0}^{(m-i)}=
\sum_{\bm{c}\in\mathbb{N}^{\ell_1}}\sum_{\lambda\neq\lambda_0,\ \mu\neq\mu_0}
x^{\bm{c}}g_{\lambda,\mu,\bm{c}}\eta_{\mu}^{(m-i)}.
\end{align*}
for some $g_{0,\bm{c}}\in\mathbb{K}$
and for some $g_{\lambda,\mu,\bm{c}}\in S_2$.
There exists $\bm{c}\in\mathbb{N}^{\ell_1}$
such that $g_{0,\bm{c}}\neq 0$, while
\begin{align*}
\eta_{\mu_0}^{(m-i)}=\sum_{\lambda\neq\lambda_0,\ \mu\neq\mu_0}
\left(g_{\lambda,\mu,\bm{c}}/g_{0,\bm{c}}\right)\eta_{\mu}^{(m-i)},
\end{align*}
i.e., $\eta_{\mu_0}^{(m-i)}$ is a linear combination of
$\{\eta_{\mu}^{(m-i)}\}_{\mu\neq \mu_0}$ over $S_2$.
This contradicts to the minimality of $\{\eta_{\mu}^{(m-i)}\}_{\mu\in M}$.
Hence $B_i$ is a minimal set of homogeneous generators for
$D^{(i)}(\mathscr{A}_1)D^{(m-i)}(\mathscr{A}_2)$ over $S$.

By Proposition \ref{prop-decomp-D(A1timesA2)}, we also have that
$\bigcup_{i=0}^{m}B_i$ is a minimal set of homogeneous generators for
$D^{(m)}(\mathscr{A}_1\times\mathscr{A}_2)$ over $S$.
By \cite[Theorem A.19]{Orlik-Terao},
a minimal set of homogeneous generators for a free graded module is a basis.
Thus the assumption
(i.e., $D^{(m)}(\mathscr{A}_1\times\mathscr{A}_2)$ is free over $S$)
implies that $\bigcup_{i=0}^{m}B_i$ is a basis
for $D^{(m)}(\mathscr{A}_1\times\mathscr{A}_2)$ over $S$.
Here we prove that $\{\eta_{\mu}^{(m-i)}\}_{\mu\in M}$
is linearly independent over $S_2$ for any $i$. Indeed, we assume that
\begin{align}\label{eq-lin-indepS2}
\sum_{\mu\in M}g_{\mu}\eta_{\mu}^{(m-i)}=0
\end{align}
for some $g_{\mu}\in S_2$. We fix $\lambda\in\Lambda$, and
by multiplying $\theta_{\lambda}^{(i)}$ to the equation \eqref{eq-lin-indepS2},
\begin{align*}
\sum_{\mu\in M}g_{\mu}\theta_{\lambda}^{(i)}\eta_{\mu}^{(m-i)}=0.
\end{align*}
Since $\{\theta_{\lambda}^{(i)}\eta_{\mu}^{(m-i)}\}_{\mu\in M}$
is linearly independent over $S$, we have $g_{\mu}=0$ for any $\mu\in M$.
Therefore, $\{\eta_{\mu}^{(m-i)}\}_{\mu\in M}$ is a basis for
$D^{(m-i)}(\mathscr{A}_2)$ over $S_2$, and similarly
$\{\theta_{\lambda}^{(i)}\}_{\lambda\in\Lambda}$ is a basis for
$D^{(i)}(\mathscr{A}_1)$ over $S_1$ for $0\leq i\leq m$.
\hfill$\Box$
\vspace{4mm}

\noindent
{\it Proof of Proposition \ref{prop-A_X-free}.}
Let $Q_X=Q(\mathscr{A}_X)$ and let $Q_0=Q(\mathscr{A})/Q_X$.
Since $\mathbb{K}$ is an infinite field, we can take
$w\in X\setminus\bigcup_{H\in\mathscr{A}\setminus\mathscr{A}_X}H$.
We note that $\alpha_H(w)=0$ if and only if $H\in\mathscr{A}_X$.
In other words, $Q_0(w)\neq 0$.
We suppose that $e_1,\dots,e_{\ell}$ is a basis for $V$
dual to $x_1,\dots,x_{\ell}$, and that $w=\sum_{i=1}^{\ell}w_i e_i$.
We define a $\mathbb{K}$-algebra isomorphism $\tau:S\rightarrow S$
by $\tau(x_i)=x_i+w_i$ for $1\leq i\leq \ell$.
The inverse of $\tau$ is a $\mathbb{K}$-algebra isomorphism
$\tau^{-1}:S\rightarrow S$ given by $\tau^{-1}(x_i)=x_i-w_i$
for $1\leq i\leq \ell$.
We also define a map $\tau:D^{(m)}(S)\rightarrow D^{(m)}(S)$ by
$\tau\left(\sum_{|\bm{a}|=m}h_{\bm{a}}\partial^{\bm{a}}\right)=
\sum_{|\bm{a}|=m}\tau(h_{\bm{a}})\partial^{\bm{a}}$.
Let $\theta\in D^{(m)}(S)$ and let $f\in S$.
We note that $(\tau\theta)(\tau f)=\tau(\theta(f))$ and $\tau Q_X=Q_X$.
Moreover, if $\theta\in D^{(m)}(\mathscr{A})$, then
$(\tau\theta)(Q_X f)=(\tau\theta)(\tau Q_X \tau\tau^{-1}f)
=\tau(\theta(Q_X\tau^{-1}f))\in Q_X S$ by Proposition \ref{prop-holm2.4}.
Thus $\tau D^{(m)}(\mathscr{A})\subseteq D^{(m)}(\mathscr{A}_X)$.

Let $\{\theta_1,\dots,\theta_s\}$ be a basis for $D^{(m)}(\mathscr{A})$.
We may suppose that $\det M_m(\theta_1,\dots,\theta_s)=Q(\mathscr{A})^t$. Then
\begin{align*}
\det M_m(\tau\theta_1,\dots,\tau\theta_s)=\tau\left(Q(\mathscr{A})^t\right)
=(\tau Q_X)^t (\tau Q_0)^t=Q_X^t(\tau Q_0)^t.
\end{align*}
We write $\tau\theta_i=\sum_{k\geq 0}\phi_i^{(k)}$, where
$\phi_i^{(k)}\in D^{(m)}(\mathscr{A}_X)$ and $\deg \phi_i^{(k)}=k$.
Since $\tau Q_0(0)=Q_0(w)\neq 0$, there exist
$\phi_1^{(k_1)},\dots,\phi_{s}^{(k_{s})}$ such that
$\det M_m(\phi_1^{(k_1)},\dots,\phi_{s}^{(k_{s})})=cQ_X^t$
for some $c\in\mathbb{K}\setminus\{0\}$.
Therefore, $\phi_1^{(k_1)},\dots,\phi_{s}^{(k_{s})}$
form a basis for $D^{(m)}(\mathscr{A}_X)$ by Saito's criterion.
\hfill$\Box$
\vspace{4mm}

\noindent
{\it Proof of Proposition \ref{nega-ans-Q1-irr}.}
We recall that $Q(\mathscr{A})=xyzw(x+y+z)(x+y+z+w)$.
We set $X=\{x=0\}\cap\{y=0\}\cap\{z=0\}\in L(\mathscr{A})$.
Then $Q(\mathscr{A}_X)=xyz(x+y+z)$.
We can describe that
\begin{align*}
(\mathscr{A}_X,V)=(\mathscr{A}_X,X)\times(\Phi_1,V/X).
\end{align*}
Since $(\mathscr{A}_X,X)$ is generic,
$(\mathscr{A}_X,X)$ is not $1$-free (see \cite{Orlik-Terao,Yuz}).
Then Theorem \ref{thm-prod-arr} implies that
$(\mathscr{A}_X,V)$ is not $m$-free for $m\geq 1$.
Therefore, by the contraposition of Proposition \ref{prop-A_X-free},
$(\mathscr{A},V)$ is not $m$-free for $m\geq 1$.
\hfill$\Box$

\subsection{Proof of Theorem \ref{thm-shi-not2free}}
Let us recall that
$D^{(m)}(\mathscr{A})_i=\{\theta\in D^{(m)}(\mathscr{A})\mid \deg(\theta)=i\}$
for $i\geq 0$.
\begin{Lem}\label{lem-deg=1,2}
Let $\mathscr{A}=\{H_1,\dots,H_n\}$ be an $\ell$-arrangement
with $H_i=\{x_i=0\}$ $(1\leq i\leq \ell)$.
\begin{enumerate}
\item[(1)]
Let $m\geq 1$. Then $D^{(m)}(\mathscr{A})_0=\{0\}$.
\item[(2)]
Let $m\geq 2$. Then $D^{(m)}(\mathscr{A})_1=\{0\}$
if and only if for any $1\leq i\leq \ell$, there exists
$H\in\mathscr{A}\setminus\{H_1,\dots,H_{\ell}\}$ ($H$ depends on $i$) such that
the coefficient of $x_i$ in $\alpha_H$ is not zero.
\end{enumerate}
\end{Lem}
\noindent
{\it Proof.}
(1)\ Let $\theta=\sum_{\bm{a}\in\mathbb{N}^{\ell},|\bm{a}|=m}
\lambda_{\bm{a}}\partial^{\bm{a}}\in D^{(m)}(\mathscr{A})_0$,
where $\lambda_{\bm{a}}\in\mathbb{K}$.
For any $\bm{a}=(a_1,\dots,a_{\ell})\in\mathbb{N}^{\ell}$ with $|\bm{a}|=m$,
there exists an index $i$ such that $a_i\neq 0$. Then
\begin{align*}
\bm{a}!\lambda_{\bm{a}}=\theta(x^{\bm{a}})
=\theta(x_i\cdot x^{\bm{a}-\bm{e}_i})\in x_i S.
\end{align*}
This implies that $\lambda_{\bm{a}}=0$
for any $\bm{a}=(a_1,\dots,a_{\ell})\in\mathbb{N}^{\ell}$ with $|\bm{a}|=m$.
Hence $\theta=0$.

(2)\ If there exists $i$ $(1\leq i\leq\ell)$ such that
for any $H\in\mathscr{A}\setminus\{H_1,\dots,H_{\ell}\}$,
the coefficient of $x_i$ is zero in $\alpha_H$,
then $x_{i}\partial_{i}^m\in D^{(m)}(\mathscr{A})_1$.
Hence $D^{(m)}(\mathscr{A})_1\neq \{0\}$.

Conversely let $\theta=\sum_{\bm{a}\in\mathbb{N}^{\ell},|\bm{a}|=m}
f_{\bm{a}}\partial^{\bm{a}}\in D^{(m)}(\mathscr{A})_1$,
where $\deg(f_{\bm{a}})=1$.
Let $\bm{a}=(a_1,\dots,a_{\ell})\in\mathbb{N}^{\ell}$ with $|\bm{a}|=m$.
We first assume that at least two entries are nonzero in $\bm{a}$.
Let $a_i,a_j$ be such entries. Since
$\bm{a}!f_{\bm{a}}=\theta(x_i\cdot x^{\bm{a}-\bm{e}_i})\in x_i S$
and $\bm{a}!f_{\bm{a}}=\theta(x_j\cdot x^{\bm{a}-\bm{e}_j})\in x_j S$,
we have $f_{\bm{a}}\in x_i x_j S$. Hence $f_{\bm{a}}=0$.
Next, we assume that only one entry is nonzero in $\bm{a}$.
Let $a_i$ be the nonzero entry.
Since $\bm{a}!f_{\bm{a}}=\theta(x_i\cdot x^{\bm{a}-\bm{e}_i})\in x_i S$,
we have $f_{\bm{a}}\in x_i \mathbb{K}$.
By the arguments above, we can describe that
$\theta=\lambda_1 x_1\partial_1^m+\cdots+
\lambda_{\ell} x_{\ell}\partial_{\ell}^m$
for some $\lambda_1\dots,\lambda_{\ell}\in\mathbb{K}$.
By the assumption, for any $1\leq i\leq \ell$, there exists
$H\in\{H_1,\dots,H_{\ell}\}$ ($H$ depends on $i$) such that
the coefficient of $x_i$ in $\alpha_H$ is not zero. Then
\begin{align*}
m!x_i\lambda_i=\theta(\alpha_H x_i^{m-1})\in \alpha_H S,
\end{align*}
and hence $\lambda_i=0$.
We obtain $\theta=0$ as required.
\hfill$\Box$

\begin{Examp}
{\rm
(1)\ Let $\mathscr{A}$ be a $3$-arrangement defined by
$Q(\mathscr{A})=xyz(x-y)(x-z)(y-z)(x-y-z)$.
Since the coefficients of $x$, $y$ and $z$ in $\alpha_H=x-y-z$ are not zero.
By Lemma \ref{lem-deg=1,2}, $D^{(m)}(\mathscr{A})_1=\{0\}$ for $m\geq 2$.

(2)\ Let $\mathscr{A}$ be a $3$-arrangement defined by
$Q(\mathscr{A})=xyz(x-y)$.
Since $z\partial_z^m\in D^{(m)}(\mathscr{A})$,
$D^{(m)}(\mathscr{A})_1\neq\{0\}$ for $m\geq 2$.
}
\end{Examp}

\begin{Cor}\label{lem-deg=1,2-irr}
Let $\ell\geq 2$ and $m\geq 2$.
If $\mathscr{A}$ is irreducible, then $D^{(m)}(\mathscr{A})_i=\{0\}$
for $i=0,1$.
\hfill$\Box$
\end{Cor}

We now recall that ${\rm Shi}_{\ell}$ is the arrangement defined by
the equation \eqref{eq-shi-arr-A}.
\begin{Examp}
{\rm
\begin{align*}
Q({\rm Shi}_2)&=zx_1 x_2(x_1-z)(x_2-z)(x_1-x_2)(x_1-x_2-z),\\
Q({\rm Shi}_3)&=
\begin{array}{l}
zx_1 x_2 x_3(x_1-z)(x_2-z)(x_3-z)(x_1-x_2)(x_1-x_3)(x_2-x_3)\\
\times(x_1-x_2-z)(x_1-x_3-z)(x_2-x_3-z).
\end{array}
\qquad\Box
\end{align*}
}
\end{Examp}
We set $X=\{x_1=0\}\cap\{x_2=0\}\cap\{z=0\}\in L({\rm Shi}_{\ell})$.
Then
\begin{align}
Q\left(({\rm Shi}_{\ell})_{X}\right)=
zx_1 x_2(x_1-z)(x_2-z)(x_1-x_2)(x_1-x_2-z).
\end{align}
In other words, $\left({\rm Shi}_{\ell}\right)_{X}={\rm Shi}_2\times\Phi_{\ell-2}$.
By Theorem \ref{thm-prod-arr} and Proposition \ref{prop-A_X-free},
the following claim is enough to verify Theorem \ref{thm-shi-not2free}.

\begin{Claim}\label{claim-shi2-not-2free}
The arrangement ${\rm Shi}_2$ is not $2$-free.
\end{Claim}

\noindent
{\it Proof of Claim \ref{claim-shi2-not-2free}.}
We describe that $Q({\rm Shi}_2)=xyz(x-y)(x-z)(y-z)(x-y-z)$.
Let
\begin{align*}
\theta_E=&\,x^2\partial_x^2+y^2\partial_y^2+z^2\partial_z^2+
2xy\partial_x \partial_y+2xz\partial_x \partial_z+2yz\partial_y \partial_z,\\
\theta_1=&\,x(x-z)(x-y)(x-y-z)\partial_x^2,\\
\theta_2=&\,y(x-y)(y-z)(x-y-z)\partial_y^2,\\
\theta_3=&\,z(x-z)(y-z)(x-y-z)\partial_z^2,\\
\theta_4=&\,xy(x-z)(y-z)\left(\partial_x+\partial_y\right)^2,\\
\theta_5=&\,
\begin{array}{l}
y^2(x-y)(y-z)\partial_y^2-z^2(x-z)(y-z)\partial_z^2
+xy(x-y)(y-z)\partial_x \partial_y\\
-xz(x-z)(y-z)\partial_x \partial_z-yz(y-z)^2\partial_y \partial_z.
\end{array}
\end{align*}
Then we can directly check that
$\theta_E,\theta_1,\dots,\theta_5\in D^{(2)}({\rm Shi}_2)$,
using Corollary \ref{cor-check-D(A)}. Since
\begin{align*}
\det M_m(\theta_E,\theta_1,\dots,\theta_5)=4(y-z)Q^3\neq 0,
\end{align*}
the operators $\theta_E,\theta_1,\dots,\theta_5$
are linearly independent over $S$.
Meanwhile, by Lemma \ref{lem-deg=1,2} and by solving linear equations,
the vector spaces $D^{(2)}({\rm Shi}_2)_i\ (0\leq i\leq 3)$ are determined
as follows:
\begin{align*}
&D^{(2)}({\rm Shi}_2)_0=\{0\},\\
&D^{(2)}({\rm Shi}_2)_1=\{0\},\\
&D^{(2)}({\rm Shi}_2)_2=\mathbb{K}\theta_E,\\
&D^{(2)}({\rm Shi}_2)_3=
(\mathbb{K}x+\mathbb{K}y+\mathbb{K}z)\theta_E.
\end{align*}
Thus $\theta_1,\dots,\theta_5$ cannot be expressed as
an $S$-linear combination of operators in $D^{(2)}({\rm Shi}_2)$
of degree $\leq 3$.
Therefore, by the linear independence for $\theta_E,\theta_1,\dots,\theta_5$
over $S$, there exists a minimal set $G$ of generators
for $D^{(2)}({\rm Shi}_2)$ such that
$\{\theta_E,\theta_1,\dots,\theta_5\}$ is contained in $G$.

We now assume that ${\rm Shi}_2$ is $2$-free.
Then the rank of $D^{(2)}({\rm Shi}_2)$ is equal to six, and hence
$G=\{\theta_E,\theta_1,\dots,\theta_5\}$ is a basis for $D^{(2)}({\rm Shi}_2)$.
However, the fact that
$\det M_m(\theta_E,\theta_1,\dots,\theta_5)=4(y-z)Q^3
\not\in\mathbb{K}Q^3\setminus\{0\}$ contradicts to Saito's criterion.
\hfill$\Box$

\section{Open problems}
We conclude this paper by giving three open questions.
Theorem \ref{thm-shi-not2free} states that
${\rm Shi}_{\ell}$ is not $2$-free for any $\ell\geq 2$.
However, we can show that ${\rm Shi}_{2}$ is $m$-free for $3\leq m\leq 7$
by a straightforward calculation, using MATLAB.
It seems to be interesting to determine
$m$-freeness of ${\rm Shi}_{\ell}$.
\begin{Prob}
{\rm
Determine $m$-freeness of ${\rm Shi}_{\ell}$ for all $m\geq 3$
and all $\ell\geq 2$.
}
\end{Prob}
There are two other interesting problems which have not answered yet.
\begin{Prob}
{\rm
Are supersolvable arrangements $m$-free for all $m\geq 1$?
}
\end{Prob}
\begin{Prob}
{\rm
Give a sufficient condition for arrangements to be $m$-free
for all $m\geq 1$.
}
\end{Prob}

\section*{Acknowledgments}
The first author is partially supported by
JSPS Grant-in-Aid for Scientific Research (B)
16H03924, and Grant-in-Aid for Exploratory Research 16K13744.
The second author is supported by
JSPS Grant-in-Aid for Young Scientists (B) 16K17582.

\end{document}